\documentclass [12pt] {article}
\usepackage{amsfonts}
\newcommand{\qed} {\hspace {0.1in} \rule {1.5mm} {3.5mm}}

\newtheorem{lemma}{Lemma}[section]
\newtheorem{corollary}{Corollary}[section]
\newtheorem{theorem}{Theorem}
\newtheorem{proposition}{Proposition}[section]
\newtheorem{definition}{Definition}[section]

\def\e{\epsilon}

\def\G{\Gamma}

\def\ov{\overline{V}}
\def\r{\mbox{rank}\,}

\def\sq{\subseteq}

\def\dim{{\rm dim}}

\def\Fo{F$\mbox{\o}$lner}
\def\<{\langle}
\def\>{\rangle}
\def\proof{\smallskip\noindent{\bf Proof:} }

\def\bR{{\mathbb R}}
\def\bN{{\mathbb N}}

\def\cf{\mbox{$\cal F$}}

\def\to{\rightarrow}

\title{The amenability of affine algebras}
\author{{\sc G\'abor Elek}
\cr Mathematical Institute of
the Hungarian Academy of Sciences\cr P.O. Box 127, H-1364 Budapest, Hungary\cr
elek@renyi.hu}
\date{}
\begin{document}

\maketitle
\noindent{\bf Abstract.}  We introduce the notion of amenability for affine
algebras. We characterize amenability by \Fo-sequences, paradoxicality and
the existence of finitely invariant dimension-measures. Then we extend the
results of Rowen on ranks, from affine algebras of
subexponential growth to amenable
affine algebras. 
\vskip 0.2in
\newpage
\section{Introduction}
First, let us recall the classical notion of amenability. Let $\Gamma$ be a
discrete group. We call $\G$  paradoxical
, if it can be written as a disjoint
union $\G=A_1\cup A_2\cup\dots\cup A_m$ such that for some elements
$g_1,h_1,g_2,h_2,\dots, g_m,h_m\in\G,$ the sets $A_1g_1, A_1h_1, A_2g_2,
A_2h_2,\dots,A_m g_m,$

\noindent
 $A_mh_m$ are disjoint as well. The group $\G$ is called
 amenable if it is not paradoxical. The theorem below is one the
fundamental results on amenability.
\begin{theorem}
The following conditions are equivalent.
\begin{enumerate}
\item $\G$ is amenable.
\item There exists a finitely additive measure on the subsets of $\G$ such
that $\mu(\G)=1$ and $\mu(Ag)=\mu(A)$ for any $A\sq\G$ and $g\in\G$.
\item There exists a sequence of finite subsets (\Fo-exhaustion)

\noindent
 $\cf_1\sq\cf_2\sq\dots,\cup_{n=1}^\infty \cf_n=\G$, such that for any
 $g\in\G$,
$$\lim_{n\to\infty}\frac{|\cf_n\cup\cf_n g|}{|\cf_n|}=1\,.$$ \end{enumerate}
\end{theorem}
 The goal of this
paper is to define and study the appropriate version of amenability for
affine algebras. Throughout this article $R$ denotes an affine algebra (not
necessarily unital) over a commutative field $K$.
\begin{definition}
\label{d2}
The affine algebra $R$ is (left) amenable if there exists a sequence of finite
dimensional linear
subspaces $W_1\sq W_2\sq\dots,\cup_{n=1}^\infty W_n=R$, such that for
any $r\in R$
\begin{equation}
\label{e2}
\lim_{n\to\infty} \frac{\dim_K(W_n r+ W_n)} {\dim_K (W_n)}=1\,.
\end{equation}

\end{definition}
We call such an exhaustion by subspaces a \Fo-exhaustion.
Now we define the analogues of paradoxicality and the invariant finitely
additive measure for algebras without zero-divisors.
\begin{definition}
Let $R$ be an affine algebra without zero divisors. We say that $R$ is
paradoxical, if {\it any} basis of $R$ over $K$, $\{f_i\}_{i=1}^\infty$ can be written
as the disjoint union $A_1\cup A_2\cup\dots\cup A_m$ such that for some
non-zero elements $g_1,h_1,g_2,h_2,\dots,g_m,h_m\in R$, the sets \\
$A_1g_1,A_1h_1,A_2g_2,A_2h_2,\dots,A_mg_m,A_mh_m$ are mutually independent.
\end{definition}
Now let $\{e_i\}_{i=1}^\infty$ be a basis of $R$, where again $R$ has no
zero-divisors.
An independent subset $L\subset R$ is called regular with respect to
$\{e_i\}_{i=1}^\infty$ if there exists subsets of $\{e_i\}_{i=1}^\infty$:\,
$A_1,A_2,A_3,\dots,A_n$ and $\{r_1,r_2,\dots,r_n\}\subset R$ such that
$L$ can be written as the disjoint union of $A_1r_1,A_2r_2,\dots,A_nr_n$.
\begin{definition}
An invariant finitely additive dimension-measure with respect
to $\{e_i\}_{i=1}^\infty$  is a non-negative function $\mu$ on the set of regular
subsets satisfying the following conditions:
\begin{enumerate}
\item 
$\mu(\{e_i\}_{i=1}^\infty)=1$ and $\mu(A)\leq 1$ for any independent regular
subset $A$.
\item
If $A$ and $B$ are independent regular subsets then
$\mu(A\cup B)=\mu(A)+\mu(B)$.
\item
For any non-zero $r\in R$ and regular set $A$, $\mu(A)=\mu(Ar)$.
\end{enumerate}
\end{definition}
 The main result of the paper is that the following theorem.
\begin{theorem}
The following conditions are equivalent for affine algebras $R$ without
zero-divisors.
\begin{enumerate}
\item  $R$ is amenable.
\item $R$ is not paradoxical.
\item  There exists a finitely additive invariant dimension-measure on $R$ with
respect to some basis $\{e_i\}_{i=1}^\infty$.
\end{enumerate}
\end{theorem}
We shall also study the algebraic properties of amenable algebras, extending
Rowen's work on algebras of subexponential growth e.g. we prove that amenable
affine algebras has the unique rank property.
\section{The proof of Theorem 2.}
\subsection{The Doubling Lemma}
\begin{lemma}
\label{3.1l}
Let $R$ be a non-amenable affine algebra with no zero-divisor. Then
there exists a finite dimensional linear subspace
$Z\subset R$ and $\epsilon>0$ such that for any finite dimensional linear
subspace $V\subset R$,
$$\frac{\dim_K(VZ+V)}{\dim_K(V)}>1+\epsilon\quad.$$
\end{lemma}
\proof
Let $Z_1\subset Z_2\subset \dots,\quad\cup_{n=1}^\infty Z_n=R$ be
a sequence of finite dimensional subspaces. Suppose that the
statement of the Lemma is not true, then there exist finite dimensional
linear subspaces $V_1,V_2,\dots,V_n,\dots$ such that

$$\frac{\dim_K(V_nZ_n+V_n)}{\dim_K(V_n)}<1+\frac{1}{2^n}\quad.$$
Obviously, $\dim_K(V_n)\to\infty$. Now we construct a \Fo-exhaustion for $R$
inductively.
Let $W_1=V_1$. If we have already constructed
$W_1\subset W_2\subset\dots\subset W_{n-1}\,$ then choose
$k$ such a way that
$$\dim_K(V_k)\geq (\dim_K(W_{n-1})+\dim_K(Z_n))\cdot 2^n\,.$$
Let $W_n=V_k+W_{n-1}+Z_n\,.$ Then $\{W_n\}_{n=1}^\infty$ will satisfy
(\ref{e2}), leading to a contradiction.\quad\qed\vskip 0.2in

\noindent
As a corollary we have the following Doubling Lemma.
\begin{lemma}
\label{3.2l}
Let $R$ be a non-amenable affine algebra with no zero-divisor. Then
there exists a finite dimensional linear subspace
$Z\subset R$ such that for any finite dimensional linear
subspace $V\subset R$,
$$\frac{\dim_K(VZ)}{\dim_K(V)}>2\quad.$$
\end{lemma}

\subsection{Amenability implies the existence of finitely additive invariant
dimension measure}
\begin{lemma}\label{l50}
Let $R$ be an amenable affine algebra with no zero divisor.
Then one can construct a sequence of finite dimensional vector spaces,
$\ov_1\subset V_1\subset \ov_2\subset V_2\subset\dots\subset R $
with the following properties.
\begin{itemize}
\item $\{V_n\}_{n=1}^\infty$ satisfy (\ref{e2}).
\item $\lim_{n\to\infty}\frac{\dim_K(\ov_n)}{\dim_K(V_n)}=1$\,.
\item
For any finite dimensional linear subspace $Z\subset R$
there exists $k>0$, such that if $n>k$, then $\ov_n+\ov_nZ\subset V_n$.
\end{itemize}
\end{lemma}
\proof
First choose an exhaustion
$W_1\subset W_2\subset\dots,\,\cup_{n=1}^\infty W_n=R$ satisfying (\ref{e2}).
Let $\ov_1=W_1,\,V_1=W_1$. Suppose that we have already chosen
$$\ov_1\subset V_1\subset \ov_2\subset V_2\dots\subset
\ov_{n-1}\subset V_{n-1}\,.$$
Let us pick $k$ so large that $V_{n-1}\subset W_k$. Then we choose $l>k$
so that
$$\frac{\dim_K((W_l W_k+W_l)W_k+(W_l W_k+W_l))}{\dim_K(W_l)}\leq1+\frac
{1}{2^n}\,.$$
Then let $\ov_n=W_l$ and $V_n=W_l W_k+W_l\,.$ \quad\qed
\begin{proposition}
\label{p51}
Let $R$ be an amenable affine algebra with no zero divisor. Then there exists
a finitely additive invariant dimension-measure on $R$ with respect to
some basis $\{e_i\}^\infty_{i=1}\,.$
\end{proposition}
\proof
Let us choose a basis $\{e_i\}^\infty_{i=1}$ of $R$ inductively, such a way
that if $\ov_i$ is a $k_i$-dimensional space, then $\{e_1,e_2,\dots,e_{k_i}\}$
form a basis of $\ov_i$, similarly if
$V_i$ is a $l_i$-dimensional space, then $\{e_1,e_2,\dots,e_{l_i}\}$
form a basis of $V_i$.
\begin{lemma}
\label{vl}
For $0\neq s\in R$, let
$$F_k(s)=\{e_j\in V_k\,:\,e_js\notin V_k\}$$
$$B_k(s)=\{e_j\notin V_k\,:\,e_js\in V_k\}$$
Then,
\begin{equation}
\label{edo1}
\lim_{k\to\infty}\frac{|F_k(s)|}{\dim_K(V_k)}=0
\end{equation}
\begin{equation}
\label{edo2}
\lim_{k\to\infty}\frac{|B_k(s)|}{\dim_K(V_k)}=0
\end{equation}
\end{lemma}
\proof
First note that if $k$ is large, then
$\ov_ks\subset V_k$, consequently (\ref{edo1}) holds. Then (\ref{edo2})
follows from the fact that the right multiplication by $s$ is an injective
map.
\qed

\vskip 0.2in
Now we define the finitely additive invariant dimension-measure. For any
regular independent subset $L$,
$\mu(L)=\lim_{\omega}\frac{|L\cap V_k|}{\dim_K(V_k)}\,.$ Then of course,
$\mu(L)\leq 1$,
$\mu(\{e_i\}_{i=1}^\infty)=1$ and $\mu(A)+\mu(B)=\mu(A\cup B)$ if
$A$ and $B$ are independent. In order to finish the proof of
Proposition \ref{p51}
it is enough to see that for any $0\neq r\in R$ and regular independent
subset $L$,
\begin{equation}
\label{edo3}
\lim_{k\to\infty} \frac{|\dim_K(Lr\cap V_k)-\dim_K(L\cap V_k)|}
{\dim_K(V_k)}=0\,.
\end{equation}
However, by additivity, we may suppose that $L$ is constructed by using 
only one translation, that is for any $a_i\in L$ there exists
$e_{n_i}$ such that $a_i=e_{n_i} s\,.$
Let $N_L\subset\{e_i\}^\infty_{i=1}$ be the
set of all such $e_{n_i}$'s. Then $L=N_L s$. By Lemma \ref{vl},
$\mu(N_L)=\mu(N_Ls)$ and $\mu(N_L)=\mu(N_Lsr)$ that implies the
invariance of $\mu$. \quad\qed

\subsection{Non-amenability implies paradoxicality}
The goal of this subsection is to prove the following proposition.
\begin{proposition}
\label{p23}
If the amenable affine algebra is not amenable then it is paradoxical.
\end{proposition}
We apply the ``algebraization'' of the tools used in \cite{DSS}.
Our first lemma is just the linear algebraic analog of the classical
Hall lemma of graph theory.
\begin{lemma}\label{l23}
Let $e_1,e_2,\dots,e_m$ be a basis for the $m$-dimensional vector space $K^m$ 
and let $T_1,T_2,\dots,T_k$ be a finite collection of linear transformations 
from $K^m$ to $K^n$.
Suppose that for any $l$-tuple $\{e_{i_1},e_{i_2},\dots,e_{i_l}\}$, the linear
vector space spanned by the vectors \\
 $\{\cup_{t=1}^l\cup_{j=1}^k T_j(e_{i_t})\} $
is at
least $l$-dimensional.
Then, there exists a \\ function $\phi:\{1,2,\dots,m\}\to
\{1,2,\dots,k\}$ such that the vectors \\ $\{T_{\phi(1)}(e_1),T_{\phi(2)}(e_2),
\dots T_{\phi(m)}(e_m)\}$ are independent.
\end{lemma}
\proof 
We proceed by induction. The lemma obviously holds for $m=1$. Suppose that the
lemma holds for any $1\leq k <m$.

\noindent
If for any $l$-tuple $l<m$, $\{e_{i_1},e_{i_2},\dots,e_{i_l}\}$, the linear
vector space spanned by the vectors $\{\cup_{t=1}^l\cup_{j=1}^k
T_j(e_{i_t})\} $
is at
least $l+1$-dimensional, then first define $\phi(1)$ such a way that
$T_{\phi(1)}(e_1)$ is non-zero. Then for the remaining basis vectors 
$\{e_2,e_3,\dots,e_m\}$ let us consider the quotient maps 
$T'_j:K^{m-1}\to K^n/\{T_{\phi(1)}(e_1)\}$. This new system of vector spaces
 and
maps must satisfy the conditions of our lemma. Hence we can extend $\phi$ to
the whole set $\{1,2,\dots,m\}$.

\noindent
Now, if for some $l$-tuple $\{i_1,i_2,\dots,i_l\}$ $l<m$, the linear
vector space spanned by the vectors $\{\cup_{t=1}^l\cup_{j=1}^k
T_j(e_{i_t})\} $
is exactly $l$-dimensional, then first define $\phi$ for
 $\{i_1,i_2,\dots,i_l\}$. Then for the remaining vectors, we can again consider
the quotient maps \\
$T_j':K^{m-l}\to K^n/\{T_{\phi(i_1)}(e_{i_1}), T_{\phi(i_2)}(e_{i_2}),\dots
T_{\phi(i_l)}(e_{i_l})\}$. Again, the new system of vector spaces and maps must
satisfy the conditions of our lemma, hence we can extend $\phi$ onto the whole
set $\{1,2,\dots,m\}$. \qed

\vskip 0.2in
\noindent
Now we have the following corollary.
\begin{lemma}\label{l24}
Let $e_1,e_2,\dots,e_m$ be a basis for the $m$-dimensional vector space $K^m$ 
and let $T_1,T_2,\dots,T_k$ be a finite collection of linear transformations 
from $K^m$ to $K^n$.
Suppose that for any $l$-tuple $\{e_{i_1},e_{i_2},\dots,e_{i_l}\}$, the linear
vector space spanned by the vectors \\
 $\{\cup_{t=1}^l\cup_{j=1}^k T_j(e_{i_t})\} $
is at
least $2l$-dimensional.
Then, there exist two functions \\
 $\phi:\{1,2,\dots,m\}\to
\{1,2,\dots,k\}$ and $\psi:\{1,2,\dots,m\}\to
\{1,2,\dots,k\}$
  such that the vectors $\{T_{\phi(1)}(e_1),T_{\phi(2)}(e_2),
\dots, T_{\phi(m)}(e_m), T_{\psi(1)}(e_1),T_{\psi(2)}(e_2),
\dots T_{\psi(m)}(e_m)\}$ are independent.
\end{lemma}
\proof 
First define $\phi$ by our previous lemma then apply the same lemma
for the quotient map
 $T'_j:K^m\to K^n/[\{T_{\phi(i_1)}(e_{1}), T_{\phi(i_2)}(e_{2}),\dots,
T_{\phi(i_n)}(e_{n})\}$
The next proposition is a simple corollary of the previous lemma and the
classical K\"onig-lemma (or compactness) argument (see also \cite{DSS}).
\begin{proposition}
\label{p25}
Let $\{e_i\}_{i=1}^\infty$ be a basis for the infinite dimensional
affine algebra $R$. Let $S=\{r_1,r_2,\dots,r_s\}$ be a set of elements
in $R$. Suppose that for any $l$-tuple $\{e_{i_1}, e_{i_2},\dots, e_{i_l}\}$
the linear vector space spanned by the vectors
$\{\cup^l_{t=1}\cup^s_{j=1}\,e_{i_t}\cdot r_j\,\}$ is at least
$2l$-dimensional. Then one has a partition of
$\{e_i\}_{i=1}^\infty=A_1\cup A_2\dots\cup A_m$ and elements \\
$g_1,h_1,g_2,h_2,\dots,g_m,h_m\in S$ such that the sets
$A_1g_1, A_1h_1, A_2g_1, A_2h_1,\dots, A_mg_m,A_mh_m$ are mutually
independent.\qed
\end{proposition}
Now we prove Proposition \ref{p23}. If $R$ is non-amenable, then by
Lemma \ref{3.2l}, for any basis $\{e_i\}_{i=1}^\infty$, there exist a
subset $\{r_1,r_2,\dots,r_s\}\subset R$ satisfying the conditions
of Proposition \ref{p25}. Consequently, $R$ is paradoxical.\quad\qed

\subsection{Paradoxicality implies the non-existence of finitely additive
invariant dimension- measure}
\begin{proposition}
\label{puj}
If $R$ is a paradoxical amenable algebra, then there is no finitely additive
dimension-measure on $R$.
\end{proposition}
\proof
Suppose that $\mu$ is a finitely additive invariant dimension-measure
with respect to the basis $\{e_i\}_{i=1}^\infty$. Then consider
the paradoxical decomposition $\{e_i\}_{i=1}^\infty=A_1\cup A_2\cup\dots
\cup A_m$ as in the definition of paradoxicality. Then 
$B=A_1g_1\cup A_1h_1\cup\dots\cup A_mg_m\cup A_mh_m$ is a regular independent
subset of dimension $2$. This is a contradiction. \quad\qed

\vskip 0.2in
\noindent
Now Theorem $2$ follows from Propositions \ref{p51}, \ref{p23} and \ref{puj}.

\section{The algebraic properties of amenable algebras}
\subsection{The basic properties}
In this section we prove some of the basic algebraic properties of the
amenable algebras.
\begin{proposition}
\label{p4}
Any affine algebra of subexponential growth is amenable
\end{proposition}
\proof
Suppose that $S=\{r_1,r_2,\dots,r_k\}\sq R$ is a generator system for $R$ that
is $R=K(r_1,r_2,\dots,r_k)$. We denote by $R_m$ the $m$-ball with respect to 
$S$ that is $R_m=\sum_{j=1}^mKS^j$.
Let $d_m=\dim_K(R_m)$. Since $R$ has subexponential growth, for any $\epsilon
>0$ there exists $C_{\epsilon}>0$ such that $d_m\leq C_\epsilon(1+\epsilon)^m$
for all $m\geq 1$.
Therefore there exists a subsequence $\{d_{m_n}\}_{n=1}^\infty$ such that
$$d_{m_n+n}\leq d_{m_n}(1+\frac{1}{2^n})\,.$$
Consequently, if $W_n=R_{m_n}$, then
$$\frac{\dim_K(W_n r+ W_n)} {\dim_K (W_n)}\leq 1+\frac{1}{2^n},\,$$
provided that $r\in\sum^n_{j=1}KS^j$\,.\qed

\noindent
On the other hand, there are amenable algebras of exponential growth. It is
easy to check that if $\G$ is a finitely generated amenable group, then the
group algebra $K\G$ is amenable. Indeed, $W_n$ can be chosen as the linear
subspace spanned by the elements of $\cf_n$, where 
$\cf_1\sq \cf_2\sq\dots$ is a \Fo-exhaustion. If
$r=k_1g_1+k_2g_2+\dots +k_sg_s\in K\G$ and $\e>0$, then for sufficiently large
$n$,
$$\frac{\dim_K(W_n r+ W_n)} {\dim_K (W_n)}\leq
\frac{\dim_K(W_n+W_ng_1 +W_ng_2 +\dots +W_ng_s)} {\dim_K (W_n)}\leq 1+\e\,.$$
As it is well-known, there are amenable groups of exponential growth. In this
case $K\G$ has exponential growth.
\begin{proposition}
\label{p6}
If $R$ is an amenable affine algebra and $R$ has no zero-divisors, then $R$
has Goldie dimension $1$, that is $R$ does not contain two independent left
ideals.
\end{proposition}
\proof
Let $I,J\triangleleft R$ be left ideals, $0\neq a\in I, 0\neq b\in J\,.$
If $n$ is large enough, then $\dim_K(W_n a\cap W_n)>\frac{1}{2}\dim_K(W_n)$
and $\dim_K(W_n b\cap W_n)>\frac{1}{2}\dim_K(W_n)$, hence $\dim_K(W_na\cap
W_nb)>0\,.$
\qed

\noindent
The previous proposition shows that the group algebra of the free group of two
generators is {\it not} amenable. Later we shall need the following technical
result
on the doubling property of non-amenable algebras.
\subsection{The ranks of finitely generated modules}
Slightly modifying the arguments of Rowen \cite{Row} we define a
real-valued rank function on finitely generated (left) modules over {\it
unital} amenable affine
algebras. Let $\omega$ be an ultrafilter and $\lim_\omega:l^\infty(\bN)\to\bR$
be the corresponding ultralimit that is a linear functional on the space of
bounded
sequences such that
$$\liminf_{n\to\infty}\{a_n\}\leq\lim_\omega\{a_n\}\leq\limsup_{n\to\infty}\{a_n\}$$
and $\lim_\omega\{a_n\}=\lim_{n\to\infty}\{a_n\}$ if
$\{a_n\}^\infty_{n=1}$ is a convergent sequence of real numbers.
Note that for any finite dimensional linear subspace $Z\sq R$
containing the unit,
$$\lim_{n\to\infty} \frac{\dim_K(W_n Z)} {\dim_K (W_n)}=1\,.$$

Let $R$ be a unital amenable affine algebra with a given sequence of subspaces
$\{W_n\}^\infty_{n=1}$ satisfying (\ref{e2}).
Suppose that $M$ is a finitely generated $R$-module such that
$M=\sum^r_{j=1} Rx_i$, where $\{x_1,x_2,\dots,x_r\}\sq M$.
Then the rank of $M$ is defined as follows,
$$\r(M)=\lim_\omega
\frac{\dim_K(W_nx_1+W_nx_2+\dots+W_nx_r)} {\dim_K (W_n)}\,.$$
We shall see that the rank function might depend on the choice of the
exhaustion $\{W_n\}^\infty_{n=1}$.
\begin{proposition}
\label{p8}
The rank defined above
 does not depend on the particular choice of the generator
system
$\{x_1,x_2,\dots,x_r\}$. Also, the rank is bounded above by the minimal number
of elements spanning $M$. 
\end{proposition}
\proof
It is enough to prove that if $Z\sq R$ is a finite dimensional linear subspace
containig the unit
then,
$$\lim_{n\to\infty}
\frac{\dim_K(\sum_{i=1}^r W_nZx_i)-\dim_K(\sum_{i=1}^r W_nx_i)
} {\dim_K (W_n)}=0\,.$$
We have the following inequalities,$$
0\leq \frac{\dim_K(\sum_{i=1}^r W_nZx_i)-\dim_K(\sum_{i=1}^r W_nx_i)
} {\dim_K (W_n)}$$
$$\leq
\sum_{i=1}^r \frac{\dim_K(W_nZx_i)-\dim_K( W_nx_i)
} {\dim_K (W_n)}
\leq
r\cdot \frac{\dim_K(W_n Z)-\dim_K( W_n)
} {\dim_K (W_n)}\,.$$
However, by amenability,
$$\lim_{n\to\infty} \frac{\dim_K(W_n Z)-\dim_K( W_n)
} {\dim_K (W_n)}=0\,\quad\qed$$ 
\begin{corollary} If $R$ is a unital amenable affine algebra then,

\noindent
\label{1c10}
\begin{enumerate}

\item
$\r(R^n)=n$, that is an amenable affine algebra always satisfies the unique
rank property. (\cite{AR},\cite{Row})
\item
If $M$ and $N$ are finitely generated $R$-modules and $M$ is either a
submodule or a homomorphic image of $N$, then 
$\r(M)\leq \r(N)$.
\item
If $N$ and $M$ are finitely generated $R$-modules, then
$\r(M\oplus N)=\r(M)+\r(N)$.
\end{enumerate}
\end{corollary}
\subsection{Exact sequences}
\begin{definition}
Let $0\to M\to N\to N/M\to 0$ be an exact sequence of finitely generated
$R$-modules and let $X=\{x_1,x_2,\dots,x_n\}$ be a system of generators for
$N$, containing a system of generators for $M$. Then the relative rank is
defined as follows:
$$\r_X(M)=\lim_\omega\frac{\dim_K(M\cap \sum_{i=1}^r
W_nx_i)}{\dim_K(W_n)}\,.$$
\end{definition}
Obviously, $\r_X(M)\geq \r(M)$.
\begin{proposition}
$\r(N)=\r(M/N)+\r_X(M)$\,.
\end{proposition}
\proof
Denote by $[x_i]$ the image of the quotient map $N\to N/M$. Then,
$$\dim_K (\sum_{j=1}^r W_nx_i)= \dim_K (\sum_{j=1}^r W_n[x_i])
+ \dim_K (M\cap \sum_{i=1}^r
W_nx_i)\,.$$
Hence the statement follows. \qed
\begin{corollary}
$\r(N)\geq\r(M/N)+\r(M)$
\end{corollary}
{\bf Example:} Let $R$ be the unital algebra generated by $1,x,y$, where
$x^2=0,xy=0$. Let $W_n$ ne the linear subspace with basis
$\{1,y,y^2,\dots,y^n,x,yx,y^2x,\dots,y^{n^2}x\}$ and let $M=Rx+Ry$, $N=R$.
Then it is easy to see that $\r(N)=1, \r(M)=0, \r(N/M)=0$.
Note however, that if the linear subspaces $W_n$ are defined as
$\{1,y,y^2,\dots,y^n,x,yx,y^2x,\dots,y^{n}x\}$, then
$\r(M)=1,\r(N/M)=0$, that is the additivity holds.
In \cite{Row} the author claims that for his rank function
\begin{equation} \label{e12}\r_S(N)\leq \r_S(M)+\r_S(N/M)\,. \end{equation}
It seems to me that there might be a gap in in his argument. The previous
example suggests that the space of the exhaustion must play a greater role,
and if (\ref{e12}) is true then an ultralimit constuction would result in an
actual additive real valued rank function on the set of finitely generated
modules over affine algebras of subexponential growth. That is
$$\r(N)=\r(M/N)+\r(M)\,.\quad\qed$$
It would  immediately
imply that 
$[R^n]=[R^m]$ in the Grothendieck group $G_0(R)$. This would be much stronger than
the unique rank property. (see \cite{Luc} for a discussion).

\end{document}